\documentclass[-aos]{imsart}
\usepackage{amssymb, epsf, amsfonts,latexsym, natbib}
\usepackage{amsmath}
\usepackage[pdftex]{color}

\newtheorem{Theorem}{Theorem}

\newtheorem{conjecture}[Theorem]{Conjecture}

\startlocaldefs

\endlocaldefs

\begin{document}

\begin{frontmatter}

% "Title of the paper"
\title{On the Hermite spline conjecture and its connection to $k$-monontone densities}
\runtitle{On the Hermite spline conjecture}

% indicate corresponding author with \corref{}
% \author{\fnms{John} \snm{Smith}\corref{}\ead[label=e1]{smith@foo.com}\thanksref{t1}}
% \thankstext{t1}{Thanks to somebody} 
% \address{line 1\\ line 2\\ printead{e1}}
% \affiliation{Some University}

\author{\fnms{Fadoua} \snm{Balabdaoui}\ead[label=e1]{fadoua@ceremade.dauphine.fr}}
\address{Centre de Recherche en Math\'ematiques de la D\'ecision, Universit\'e Paris-Dauphine,
                   Paris, France,\\ \printead{e1}}
%\ead[label=e1]{???}}
%\address{\printead{e1}}

\author{\fnms{Simon} \snm{Foucart}\ead[label=e2]{foucart@math.drexel.edu}}
\address{Drexel University, Department of Mathematics, 269 Korman Center,
3141 Chestnut Street, Philadelphia, PA 19104,\\ \printead{e2}}
\medskip

\and
\author{\fnms{Jon A.} \snm{Wellner}\thanksref{t2}\ead[label=e3]{jaw@stat.washington.edu}}
\thankstext{t2}{Supported in part by NSF Grants DMS-0804587 and 
DMS-1104832,  by  NI-AID grant 2R01 AI291968-04, and by the Alexander von Humboldt Foundation}
\address{Department of Statistics, University of Washington, Seattle, WA  98195-4322,\\ \printead{e3}}

\runauthor{Balabdaoui, Foucart, and Wellner}

\begin{abstract}

The $k$-monotone classes of densities defined on $(0, \infty)$ have been known 
in the mathematical literature but were for the first time considered from a statistical point of view by 
\cite{MR2382657} and \cite{MR2830965}. 
In these works, the authors generalized the results established for 
monotone ($k=1$)  and convex ($k=2$) densities by giving a 
characterization of the Maximum Likelihood and Least Square 
estimators (MLE and LSE) and deriving minimax bounds for rates 
of convergence.  For $k \ge 3$, the pointwise asymptotic behavior of the 
MLE and LSE 
studied by
%given by 
\cite{MR2382657} %\cite{bw10} 
would show that the MLE and LSE attain the minimax lower 
bounds in a local pointwise sense. 
%are minimax efficient
However, the theory assumes that a certain conjecture about the 
approximation error of a Hermite spline holds true. 
The main goal of the present  note is to 
show why such a conjecture cannot be true. 
We also suggest how to bypass the conjecture and rebuilt the key proofs in the limit theory of the estimators.
\end{abstract}

%\begin{keyword}[class=AMS]
%\kwd[Primary ]{}
%\kwd{}
%\kwd[; secondary ]{}
%\end{keyword}

\begin{keyword}
\kwd{conjecture}
\kwd{asymptotic distribution}
\kwd{Hermite spline}
\kwd{$k$-monotone}
\end{keyword}

\end{frontmatter}

\section{Introduction}
For an integer $k \ge 1$, a density $g_0$ defined on $(0, \infty)$ is said 
to be $k$-monotone if it is nonincreasing when $k=1$, and if 
$(-1)^j g^{(j)}_0$ is nonincreasing and convex for all $j \in \{0, \ldots, k-2 \}$ 
when $k \ge 2$.  Considering the problem of estimating a density in 
one of these classes presents several interesting features. 
As shown in 
\cite{MR2382657}  %\cite{bw07} 
and 
\cite{MR2830965}, %\cite{bw10}, 
both the MLE and LSE of a $k$-monotone density exist. 
These estimators generalize the Grenander estimator of a nonincreasing density ($k=1$) 
and the MLE and LSE of a nonincreasing and convex density ($k=2$) studied by \cite{MR1891742}. 

While it is known that the Grenander estimator in the case $k=1$ converges pointwise at the rate
$n^{1/3}$ and that the MLE of a convex nonincreasing density  converges pointwise at the rate
$n^{2/5}$,  the rate of convergence $n^{k/(2k+1)}$ for the MLE (or LSE) of a $k$-monotone density in the 
general case $k\ge3$
studied in \cite{MR2382657} depends on a key conjecture which has not yet been verified.
In fact we show here that the spline conjectures made in \cite{MR2382657}
fail to hold.  On the other hand,
\cite{MR2520591} obtained a result concerning the global rate of convergence of the 
MLE of a $k$-monotone density for a  general $k\ge3$:  they showed that the rate of 
convergence of the MLE with respect to the Hellinger metric is indeed $n^{k/(2k+1)}$.

The limit case $k =\infty$ corresponds to the intersection of all $k$-monotone classes, that is the 
class of completely monotone densities on $(0, \infty)$. 
The latter turns out to be equal to the the class of mixtures of Exponentials, 
a consequence of  Bernstein--Widder theorem, see e.g. \cite{MR2474372}. 
The nonparametric MLE of a mixture of Exponentials was considered by 
\cite{MR653523}  %\cite{jew82} 
who showed its consistency and developed an EM  algorithm to compute the estimator. 
So far, there are no results available on the limit distribution of the 
completely monotone MLE.  As noted in 
\cite{MR2382657}, one natural approach seems to study 
the behavior of the MLE in the $k$-monotone class as $k \to \infty$ and $n \to \infty$.   
Such an approach requires evidently a deep understanding of the asymptotic behavior 
of the $k$-monotone MLE and of the distance between its knots.   

For an arbitrary $k \ge 1$, the work of \cite{MR2382657} aims to give a 
general approach to derive the limit distribution of the MLE and LSE 
at a fixed point $x_0 > 0$. More precisely, their work can be seen as 
an extension of the approach used by \cite{MR1891742} in convex estimation. 
Let $g_0$ denote the true $k$-monotone density. At $x_0$, and modulo the spline conjecture, it is shown that 
\begin{equation} 
\label{Lim}
\left (
\begin{array}{c}
n^{\frac{k}{2k+1}}(\bar{g}_{n}(x_{0})-g_0(x_{0}) ) \\
n^{\frac{k-1}{2k+1}}(\bar{g}^{(1)}_{n}(x_{0})-g^{(1)}_0(x_{0}) ) \\
\vdots \\
n^{\frac{1}{2k+1}}(\bar{g}^{(k-1)}_{n}(x_{0})-g^{(k-1)}_0(x_{0}))
\end{array} \right )
\rightarrow_{d}
\left (
\begin{array}{c}
c_{0}(x_0)H^{(k)}_k(0)\\
c_{1}(x_0)H^{(k+1)}_k(0)\\
\vdots \\
c_{k-1}(x_0)H^{(2k-1)}_k(0)\\
\end{array}\right ),
\end{equation}
where $\bar{g}_{n}$ is either the MLE or LSE, and 
\begin{eqnarray*}
c_{j}(x_0) = \Bigg \{ g_0(x_0)^{k-j}
\left(\frac{(-1)^k g^{(k)}_0(x_0)}{k!}\right)^{2j+1} \Bigg \}^{\frac{1}{2k+1}},
\end{eqnarray*}
for $j=0, \ldots, k-1$. Note that the constants $c_j(x_0), \ j=0, \ldots, k-1$ 
appear also in the asymptotic minimax lower bound for $L_1$ risk 
(see \cite{MR2830965}).  Let
\begin{eqnarray*}
Y_k(t)= \left \{
\begin{array}{ll}
\displaystyle{ \int_0^t \frac{(t-s)^{k-1}}{(k-1)!}dW(s)
+ \frac{(-1)^k k!}{(2k)!}t^{2k} }, 
%\hspace{0.5cm} 
& t \geq 0\\
\displaystyle{\int_t^0 \frac{(t-s)^{k-1}}{(k-1)!}dW(s)
+ \frac{(-1)^k k!}{(2k)!}t^{2k} }, 
%\hspace{0.5cm} 
&t < 0
\end{array}
\right.
\end{eqnarray*}
where $W$ is a two-sided Brownian motion on $\mathbb R$. The process $H_{k}$ appearing in the limit \eqref{Lim} is characterized  
by the following conditions:
\begin{description}
\item (i) \hspace{0.05cm} The process $H_k$
stays everywhere above the process $Y_k$: 
\begin{eqnarray*}\label{CondHk1}
H_k(t) \geq Y_k(t), \hspace{0.5cm} t \in \mathbb{R}.
\end{eqnarray*}
\item (ii) \hspace{0.05cm} $(-1)^k H_k $ is
$2k$-convex, i.e., $(-1)^k H^{(2k-2)}_k $ exists and is convex.
\item (iii) \hspace{0.05cm} The process $H_k$ satisfies
\begin{eqnarray*}\label{CondHk2}
\int_{-\infty}^{\infty} \left(H_k(t) - Y_k(t)\right) dH^{(2k-1)}_k(t) = 0.
\end{eqnarray*}
\item (iv) \hspace{0.05cm} If $k$ is even,
$\lim_{\vert t \vert \to \infty} (H^{(2j)}_k(t)-Y^{(2j)}_k(t))=0$ for
$j=0,\ldots, (k-2)/2$; if $k$ is odd, $\lim_{t \to \infty} (H_k(t)-Y_k(t))=0$
and $\lim_{\vert t \vert \to \infty} (H^{(2j+1)}_k(t)-Y^{(2j+1)}_k(t))=0$ for $j=0, \ldots, (k-3)/2$.
\end{description}

Because there is so far no device equivalent to the switching relationship 
device used in the monotone problem (see e.g. \cite{MR822052} and also \cite{MR2829859}),  the proof by  
\cite{MR1891742} of the limit of the convex estimators is more complex and 
built in several steps. 
One of the most crucial pieces of this proof is that the 
stochastic order $n^{-1/5}$ for the distance between two knot points of the estimators in a 
small neighborhood of $x_0$. 
This result holds true under the 
assumption that the true convex density $g_0$ is twice continuously differentiable 
in a neighborhood of $x_0$ such that $g''_0(x_0) > 0$.  
 
In the monotone problem, one can also show that  the distance between 
the jump points of the Grenander estimator is stochastically bounded 
above by $n^{-1/3}$ provided that $g_0$ is continuously differentiable in a 
neighborhood of $x_0$ such that $g'_0(x_0) < 0$. These working 
assumptions can be naturally put in the following general form: the 
true $k$-monotone density is $k$-times continuously differentiable in a 
neighborhood of $x_0$ such that $(-1)^k g^{(k)}_0(x_0) > 0$.  
Thus, it seems natural that $n^{-1/(2k+1)}$ gives the general 
stochastic order for all integers $k \ge 1$.  As noted in \cite{MR2382657},
\cite{MR1429931} have, in the context of fitting a regression curve via 
splines, already conjectured that $n^{-1/(2k+1)}$ is the order of the distance 
between the knot points of their regression spline under the assumption 
that the true regression curve satisfies our same working assumptions.

In the extension of the argument of \cite{MR1891742} to an arbitrary $k$, 
we have found that there is a need to show that an envelope of a certain 
VC-class is bounded. In the next section, we describe this fact more precisely,
and give the connection to our two spline conjectures made in \cite{MR2382657}. 
In Section 3, we show that these conjecture are false for $k=3$. The argument can be 
generalized to $k \ge 4$ but the calculations rapidly become 
cumbersome. In Section 4, we give a number of suggestions for building an alternative 
proof for the limit theory of the $k$-monotone estimators.

\section{Connection to splines and the conjectures}

We begin with some notation. For  integers $m \ge 0$ and $p \ge 1$, let us denote by 
$\mathcal S_m(a_1, a_2, ..., a_p)$ the space of splines on $[a, b]$ of degree $m$ 
and internal knots $a_1< \cdots <a_p $. The points $a$ and $b$ can be seen as external knots and will be denoted by $a_0$ and $a_{p+1}$, respectively. 
% endowed by a basis $(S_1, S_2, \cdots, S_{m+p+1})$ %:= (1, x, ..., x^m, (x - \tau_1)^m_+, 
Let $f$ be a differentiable function on $[a,b]$ (differentiable on $(a,b)$ 
and to the right and left of $a$ and $b$, respectively). 
If $m = 2k-1$, $p = 2k-4$, and $a < a_1 < \cdots < a_{2k-4} < b$, 
we know that there exists a unique (Hermite) spline  
$H_k \in \mathcal S_{2k-1}(a_1,a_2, ..., a_{2k-4})$ satisfying
\begin{eqnarray*}
H_k(a_j) = f(a_j)  \ \ \textrm{and} \ \ H'_k(a_j) = f'(a_j), \ \ \ \textrm{for $j=0, \ldots, 2k-3$}. 
\end{eqnarray*}
Note that for $k=2$ the Hermite spline reduces to the cubic polynomial interpolating $f$ at $a$ and $b$. We denote by $\mathcal H_k$ the spline interpolation operator which 
assigns to $f$ its spline interpolant $H_k$. 
%In the following, 
%we use the notation $\mathcal H_k [f](x)$ to mean the value at a point $x$ of the Hermite spline interpolant of $f$.

Let $\tilde g_n$ be the LSE of the true $k$-monotone density $g_0$ 
based on $n$ i.i.d. random variables $X_1, \cdots, X_n$. 
It was shown by \cite{MR2830965} that 
$\tilde{g}_n$ exists, is unique, and is a spline of degree $k-1$. 
Let $\tilde{H}_n$ denote its $k$-fold integral, that is 
\begin{eqnarray*}
\tilde{H}_n(x) = \frac{1}{(k-1)!} \int_0^x (x-t)^{k-1} \tilde{g}_n(t) dt.
\end{eqnarray*}
The function $\tilde{H}_n$ is important due to its direct involvement 
in the characterization of the estimator $\tilde{g}_n$. More precisely, if we consider 
the $(k-1)$-fold integral of the empirical distribution $\mathbb G_n$
\begin{eqnarray*}
\mathbb Y_n(x) = \frac{1}{(k-1)!} \int_0^x (x-t)^{k-1} d\mathbb G_n(t),
\end{eqnarray*} 
then the spline $\tilde{g}_n$ of degree $k-1$ is the LSE if and only if
the following (Fenchel) conditions hold
\begin{eqnarray}\label{CharLSE}
\tilde{H}_n(x) &\ge&  \mathbb Y_n(x), \ \ \textrm{for all $x \ge 0$,} \nonumber \\
\tilde{H}_n(x) & = &  \mathbb Y_n(x), \ \ \textrm{if $x$ is knot of $\tilde g_n$}.
\end{eqnarray}   
The greater focus put on the LSE is explained by the fact  that the 
characterization in \eqref{CharLSE} is much simpler to study,   
especially when the empirical processes involved are localized 
(see \cite{MR2382657}).  However, it was shown by \cite{MR2382657} that 
understanding the asymptotics of the LSE is enough as one can use strong 
consistency of the MLE to linearize its characterization and put it in a more familiar form. 

One of the key points in the study of the asymptotics is to note 
that the characterization of the LSE implies 
$\tilde H_n(\tau) = \mathbb  Y_n(\tau)$ and $\tilde H'_n(\tau) = \mathbb Y'_n(\tau)$
for a knot $\tau$ of 
$\tilde g_n$.  
Furthermore, given $2k-2$ knots 
$\tau_0 < \cdots < \tau_{2k-3}$, $\tilde g_n$  is uniquely determined on 
$[\tau_0, \tau_{2k-3}]$ by the interpolation equalities 
$\tilde H^{(i)}_n(\tau_j) = \mathbb Y^{(i)}_n(\tau_j)$, $i =0,1$, $j=0, \ldots, 2k-3$.     
%\SF{It seems that $\tilde H'_n(\tau_0) = \mathbb Y'_n(\tau_0)$ and $\tilde H'_n(\tau_{2k-3}) = \mathbb Y'_n(\tau_{2k-3})$ are not implied by \eqref{CharLSE}?}
In other words, $\tilde H_n$ is a Hermite spline interpolant of $\mathbb  Y_n$, i.e.,
 \begin{eqnarray*}
\tilde{H}_n(x) = \mathcal H_k[\mathbb Y_n](x)
\qquad \mbox{for } x \in [\tau_0, \tau_{2k-3}].
\end{eqnarray*}
Note that in any small neighborhood of the estimation point $x_0$, 
strong consistency of the $(k-1)$-st derivative of $\tilde g_n$ combined 
with the assumption that $g^{(k)}_0(x_0) \ne 0$ guarantee that the 
number of knots in that neighborhood tends to $\infty$ almost surely 
as $n \to \infty$. Hence, finding at least $2k-2$ knots is possible with 
probability one. At this stage, we know that $\tau_{2k-3} - \tau_0 \to 0$ 
almost surely as $n \to \infty$, and our goal is to show that this 
convergence occurs with a rate equal to $n^{-1/(2k+1)}$.    
In the next section, we describe briefly the key argument in the proof of 
\cite{MR2382657}  and recall the two related spline conjectures.

\subsection{The spline conjectures}

Take an arbitrary point $\bar{\tau} \in [\tau_0, \tau_{2k-3}]$ such 
that $\bar{\tau} \notin \{\tau_0, \cdots, \tau_{2k-3} \}$. By the inequality in \eqref{CharLSE}, 
we have that 
\begin{eqnarray*}
\mathcal H_k[\mathbb Y_n](\bar{\tau}) \ge \mathbb Y_n(\bar{\tau}).
\end{eqnarray*} 
If $Y$ denotes the population counterpart of $\mathbb Y_n$, i.e., the $(k-1)$-fold integral of $g_0$ 
\begin{eqnarray*}
Y(x)= \frac{1}{(k-1)!} \int_0^x (x-t)^{k-1} g_0(t)dt,  
\end{eqnarray*}
then the latter inequality can be rewritten in the more useful form
\begin{eqnarray}\label{FondIneq}
[\mathcal H_k Y - Y](\bar{\tau}) \ge 
\mathcal H_k[Y - \mathbb Y_n](\bar{\tau}) -  [Y - \mathbb Y_n](\bar{\tau}).
\end{eqnarray} 
Both sides of the inequality can be recognized as the Hermite interpolation 
errors corresponding to the interpolated functions $Y$ and $Y - \mathbb Y_n$. 
While $Y$ is $(2k)$-times differentiable on $[\tau_0, \tau_{2k-3}]$ under our
working assumptions, the function $Y- \mathbb Y_n$ is only $(k-2)$-times 
continuously differentiable since $\mathbb Y_n$ is the $(k-1)$-st fold integral of the 
(piecewise constant) empirical distribution function $\mathbb G_n$. 
 
Taylor expansions of $Y$ and $\mathbb Y_n - Y$ up to the orders $2k$ and $k-1$, respectively, give yet 
another form for \eqref{FondIneq}. 
On $[\tau_0, \tau_{2k-3}]$, consider the functions
\begin{multline*}
f_{0}(x) =  \frac{x^{2k}}{(2k)!},
\qquad b_u(x) = \frac{(x-u)^{k-1}_+}{(k-1)!}, \quad u \in (\tau_0, \tau_{2k-3}),\\
 r(x) = \frac{1}{(2k-1)!} \int_{\bar{\tau}}^{\tau_{2k-3}}(x-t)^{2k-1}_+ (g^{(k)}_0(t) - g^{(k)}_0(\bar{\tau})) dt. \quad
\end{multline*}
%\begin{align*}
%f_{0}(x) &=  \frac{x^{2k}}{(2k)!},
%& b_u(x) &= \frac{(x-u)^{k-1}_+}{(k-1)!}, \quad u \in (\tau_0, \tau_{2k-3}),\\
%& & r(x) & = \frac{1}{(2k-1)!} \int_{\bar{\tau}}^{\tau_{2k-3}}(x-t)^{2k-1}_+ (g^{(k)}_0(t) - g^{(k)}_0(\bar{\tau})) dt. 
%\end{align*}
%$$f_{0}(x) =  \frac{x^{2k}}{(2k)!},$$
%\begin{eqnarray*}
%b_u(x) = \frac{(x-u)^{k-1}_+}{(k-1)!}, \ \ \textrm{with \ $u \in (\tau_0, \tau_{2k-3})$},
%\end{eqnarray*}
%and 
%\begin{eqnarray*}
%r(x) = \frac{1}{(2k-1)!} \int_{\bar{\tau}}^{\tau_{2k-3}}(x-t)^{2k-1}_+ (g^{(k)}_0(t) - g^{(k)}_0(\bar{\tau})) dt.  
%\end{eqnarray*}
% and its Hermite interpolant $f_t(x) = \mathcal H_k [g_t](x)$
Let $e_k =  f_0 - \mathcal H_k f_0$ be the error associated with Hermite interpolation of $f_0$.  
Then, \eqref{FondIneq} is equivalent to
\begin{eqnarray*}
g^{(k)}_0(\bar{\tau}) e_k(\bar{\tau}) \le \mathbb E_n + \mathbb R_n
\end{eqnarray*} 
where, with $G_0$ denoting the c.d.f. of $g_0$,
\begin{eqnarray*}
\mathbb E_n = \int_{\tau_0}^{\tau_{2k-3}} \mathcal H_k [b_u] (\bar{\tau}) d(\mathbb G_n(u) - G_0(u)) 
\qquad \mbox{and} \qquad
\mathbb R_n = \mathcal H_k [r](\bar{\tau}).
\end{eqnarray*}
Recalling that $(-1)^{k} g^{(k)}_0(x_0) > 0$,
so that $(-1)^k g^{(k)}_0$ is positive 
on 
a neighborhood 
$[x_0-\delta, x_0+\delta]$ for some $\delta > 0$, 
 \eqref{FondIneq} can also be rewritten as
\begin{eqnarray*}
(-1)^k g^{(k)}_0(\bar{\tau}) (-1)^k e_k(\bar{\tau}) \le \mathbb E_n + \mathbb R_n.
\end{eqnarray*}

The term $\mathbb E_n$ is an empirical 
process indexed by the class of functions $h$ such that
\begin{eqnarray*}
h(u) = h_{s, s_0, \ldots, s_{2k-3}}(u) = \mathcal H_k [b_u](s) 1_{[s_0, s_{2k-3}]} (u),
\end{eqnarray*} 
for some $s_0 < \cdots < s_{2k-3}$ in $[x_0 - \delta, x_0 + \delta ]$ and $s \in (s_0, s_{2k-3})$. 
Here $\mathcal H_k [f]$ denotes the Hermite spline interpolating 
$f$ at $s_j, j=0, \cdots, 2k-3$. The second term $\mathbb R_n$ is equal to 
the interpolation error corresponding to the $(2k)$-times differential function $r$. 
The main goals are:   (a) find upper stochastic bounds for $\mathbb E_n$ and 
$\mathbb R_n$;  (b) find a lower bound for  $(-1)^k e_k(\bar{\tau}) $ 
as a function of a power of the distance $\tau_{2k-3} - \tau_0$.  

In the absence of any knowledge about the location and distribution of the random knots $\tau_0,\ldots,\tau_{2k-3}$, 
it seems naturally desirable to get of rid of any 
dependency on these points.  
This motivates the assumption that 
the interpolation error is uniformly bounded independently of the knots. 
Thus, the following conjectures were formulated in \cite{MR2382657} to tackle (a).

\begin{conjecture}\label{Conjecture1}
Let $a =0$, $b=1$, and $b_t(x) = (x -t)^{k-1}_{+}/(k-1)!$ for $t \in (0,1)$. 
There exists a constant $d_k > 0$ such that
\begin{eqnarray}\label{Conj1}
\sup_{t \in (0,1)} \sup_{0 < y_1 < \cdots < y_{2k-4} < 1} \Vert b_t - \mathcal H_k b_t \Vert_\infty \le d_k .
\end{eqnarray}
\end{conjecture}

\begin{conjecture}\label{Conjecture2}
Let $a =0$ and $b=1$.
Then there exists a constant $c_k > 0$ such that,
for any $f \in C^{(2k)}[0,1]$,
\begin{eqnarray}\label{Conj2}
 \sup_{0 < y_1 < \cdots < y_{2k-4} < 1} \Vert f- \mathcal H_k f \Vert_\infty \le c_k \Vert f^{(2k)} \Vert_\infty.
\end{eqnarray}
\end{conjecture}

Note that Conjecture \ref{Conjecture1} cannot hold if Conjecture 
\ref{Conjecture2} does not:
indeed, in view of the Taylor expansion
$$
f(x) = \sum_{j=0}^{2k-1} f^{(j)}(0) \frac{x^j}{j!} + \int_{0}^1 f^{(2k)}(t) \frac{(x-t)_+^{2k-1}}{k!} dt
$$
and of the fact that polynomials of degree $\le 2k-1$ are preserved by  $\mathcal H_k$,
we observe that \eqref{Conj1} implies \eqref{Conj2} with $c_k = d_{2k}$.

Let us fix $s_0$ and $R > 0$ such that $[s_0, s_0 + R] \subset [x_0 - \delta, x_0 + \delta ]$. 
Conjecture~\ref{Conjecture1} implies that the class 
\begin{eqnarray*}
\mathcal F_{s_0, R} = \{ h_{s, s_0, \ldots, s_{2k-3}}: [s_0, s_{2k-3}] \subset [s_0, s_0 + R] 
\subset [x_0 - \delta, x_0 + \delta] \}
\end{eqnarray*}
admits a finite envelope, e.g.  
\begin{eqnarray*}
F_{s_0, R} (x) = a_k R^{k-1} 1_{[s_0, s_0 + R]} (x) 
\end{eqnarray*}
where $a_k > 0$ is a constant depending only on $k$ (through $d_k$). 
Together with the fact that the class $\mathcal F_{s_0, R}$ is a VC-subgraph,
this gives one 
of the most crucial results that helps establishing the stochastic order 
of the gap: the \lq\lq right\rq\rq \ stochastic bound     
\begin{eqnarray}\label{Res1}
\mathbb E_n = O_p \big( n^{-2k/(2k+1)} \big) + o_p \big( (\tau_{2k-3} - \tau_0)^{2k} \big).  
\end{eqnarray}
On the other hand, the term $\mathbb R_n$ could be bounded using Conjecture \ref{Conjecture2}. 
Since $\mathbb R_n$ is $(2k)$-times continuously differentiable on 
a neighborhood of $x_0$, Conjecture 2 yields 
\begin{eqnarray}\label{Res2}
\mathbb R_n = o_p \big( (\tau_{2k-3} - \tau_0)^{2k} \big).  
\end{eqnarray}
It follows that   
\begin{eqnarray*}
\sup_{\bar{\tau} \in [\tau_{j_0}, \tau_{j_0 +1}]} (-1)^k e_k(\bar{\tau}) 
\le O_p \big( n^{-2k/(2k+1)} \big) + o_p \big( (\tau_{2k-3} - \tau_0)^{2k} \big),
\end{eqnarray*}
where $[\tau_{j_0}, \tau_{j_0 +1}]$ is s largest knot interval 
among $[\tau_j, \tau_{j+1}], j=0, \ldots, 2k-4$. 

At this stage of the argument, the stochastic order of the gap can be 
shown to be $n^{-1/(2k+1)}$ if there exists $M > 0$ such that 
\begin{eqnarray*}
\sup_{\bar{\tau} \in [\tau_{j_0}, \tau_{j_0 +1}]} (-1)^k e_k(\bar{\tau}) > M (\tau_{2k-3} - \tau_0)^{2k}. 
\end{eqnarray*}
This can be shown using some known results on monosplines 
and Chebyshev polynomials (see \cite{Balabda-Well:05, Balabda-Well:06, MR2382657}). %Balabdaoui and Wellner (2006, 2007)).
  
Conjecturing  boundedness of the Hermite spline interpolant was a 
crucial assumption to obtain the right stochastic bound for the empirical 
process $\mathbb E_n$ and the remainder term $\mathbb R_n$. 
However, this boundedness served only as a sufficient condition. 
In the next section, we show that Conjecture 2 (hence Conjecture 1) 
is in fact answered negatively.

\section{Unboundedness of the Hermite interpolation error} 
\label{sec:InterpErrorUnbounded}

We now prove that the statement of Conjecture \ref{Conjecture2}
(and even a weaker statement where $c_k$ would be allowed to depend on $f$)
is violated for the function $f= S_*$
defined by
\begin{equation}\label{S*}
 S_*(t) = 
 S_*(t; \tau_1,\ldots,\tau_{2k-4})
 = 
 \frac{1}{(2k)!} \left ( t^{2k}
+ 2 \sum_{i=1}^{2k-4} (-1)^i (t-\tau_i)_+^{2k} \right).
\end{equation}
This choice is dictated by the fact
(not necessary here, so not proven)
that,
for $0= \tau_0 < \tau_1 < \cdots < \tau_{2k-4} < \tau_{2k-3}=1$
and for $t \in [0,1]$,
\begin{equation}
\sup_{f \in W^{2k}_\infty, \, \| f^{(2k)} \|_{\infty} \le 1} \; \Big | [\mathcal H_k f](t) - f(t) \Big|
= \Big | [\mathcal H_k S^*](t) - S^*(t) \Big| .
\label{PerfectSplineUpperBound}
\end{equation}
Setting $\mathcal E_k: = \mathcal H_k (S_*) - S_*$,
 the Landau--Kolmogorov inequality (see e.g. \cite{Kolm:39}, \cite{Kolm:62}, or \cite{MR0315070}) %Kolmogorov, 1939)
 guarantees the existence of a constant $D_{k} > 0$ depending only on $k$ such that
$$
\Vert \mathcal E^{(2k-1)}_k \Vert_\infty 
\le D_k \ \Vert \mathcal E_k \Vert_\infty^{\frac{1}{2k}} \ \Vert \mathcal E^{(2k)}_k \Vert_\infty^{\frac{2k-1}{2k}}.   
$$
Since $ \mathcal E^{(2k)}_k = - S_*^{(2k)}$ alternates between $+1$ and $-1$,
so that
$\Vert \mathcal E^{(2k)}_k \Vert_\infty  = 1$, 
it follows that if Conjecture 2 was true, then 
$\Vert \mathcal E^{(2k-1)}_k \Vert_\infty$
would be bounded independently of the knots.  
Studying the latter turns out to be easier than studying $\Vert \mathcal E_k \Vert_\infty $ itself,
 as $\mathcal E^{(2k-1)}_k$
 is a piecewise linear function (not necessarily continuous at the knots) 
 whose slope alternates between $+1$ and $-1$.  
Let us note that $\mathcal E_k$ belongs to the space 
$$
\Omega_{k}(\tau_1, \cdots, \tau_{2k-4}) = \Big \{ \gamma S_*(t) + s(t),  
\ \gamma \in \mathbb R, \ s \in \mathcal S_{2k-1}(\tau_1, \cdots, \tau_{2k-4}) \Big \},
$$  
which is a $(4k-3)$-dimensional weak Chebyshev space (see e.g. Lemma 1 in 
\cite{MR1866379}). %Bojanov and Naidenov, 2002). 
Let us also note that $\mathcal E_k$ has double zeros occurring at the 
knots $\tau_0,\tau_1, \cdots, \tau_{2k-4},\tau_{2k-3}$. 
Since $4k-4$ is the maximal number of zeros for a nonzero function in a weak Chebyshev space of dimension $4k-3$, 
there exists a constant $C \in \mathbb R$ such that,
for all $t \in [0,1]$,  
$\mathcal E_k(t)$ equals
$$
C  \left \vert
\begin{array}{ccccccc}
B_1(0) & B'_1(0)  & \cdots &  \cdots &  B_1(1) &  B'_1(1) &  B_1(t) \\ 
B_2(0) & B'_2(0)  & \cdots &  \cdots & S_2(1) & B'_2(1) &  B_2(t) \\ 
\vdots & \vdots &  \vdots & \vdots & \vdots & \vdots  &   \vdots \\
\vdots & \vdots & \vdots & \vdots & \vdots & \vdots  & \vdots \\
\vdots & \vdots &  \vdots & \vdots & \vdots & \vdots   & \vdots \\
B_{4k-3}(0) & B'_{4k-3}(0) & \cdots & \cdots  & B_{4k-3}(1) &  B'_{4k-3}(1)&  B_{4k-3}(t)
\end{array}
\right \vert ,
$$
where $(B_1, \ldots, B_{4k-3} )$ is any basis for
$\Omega_k(\tau_1, \cdots, \tau_{2k-4})$. 
The value of $C$ is determined by 
$ \mathcal E^{(2k)}_k(t) =  -1$
for $0 \le t < \tau_1$.
Our objective is now to prove the unboundedness of $\Vert \mathcal E^{(2k-1)}_k \Vert_\infty$,
which we do in the particular case $k=3$.
We consider the basis for $\Omega(\tau_1, \tau_2)$ (which has dimension 9) given by
$$
\big( 1,t,t^2,t^2(t-\tau_1), t^2(t-\tau_1)^2, t^2(t-\tau_1)^2(t-\tau_2), (t-\tau_1)^5_+, (t-\tau_2)^5_+, S_*(t) \big).  
$$ 
The determinantal expression of $\mathcal E_3(t)$ can be explicitly written as
\begin{equation}
\label{E3asDet}
\mathcal E_3(t) = C 
\begin{vmatrix}
1 & 0 & \vline & {\rm x} & \cdots & \cdots & \cdots & {\rm x} & \vline & 1\\
0 & 1 & \vline & {\rm x} & \cdots & \cdots & \cdots & {\rm x} & \vline & t\\
\hline
0 & 0 & \vline & &  &  &  & & \vline & t^2\\
0 & 0 & \vline & &  &  &  & & \vline & t^2(t-\tau_1)\\
\vdots & \vdots & \vline & &  &  & & & \vline &  t^2(t-\tau_1)^2\\
\vdots & \vdots & \vline & &  & {\Large D} & & & \vline & t^2(t-\tau_1)^2(t-\tau_2)\\
\vdots & \vdots & \vline & &  &  &  & & \vline & (t-\tau_1)^5_+\\
0 & 0 & \vline & &  &  &  & & \vline & (t-\tau_2)^5_+\\
0 & 0 & \vline & &  & &  & & \vline & S_*(t)\\
\end{vmatrix},
\end{equation}
where $D$ is 
the $7 \times 6$ matrix  
\small
%\begin{eqnarray*}  
%\left[
%\begin{array}{ccccccc}
$$\hspace{-1mm}
\begin{bmatrix}
\tau_1^2  & 2 \tau_1 & \tau_2^2   & 2 \tau_2 & 1  & 2  \\ 
0 & \tau_1^2  & \tau^2_2 (\tau_2 - \tau_1) & \tau_2 (3 \tau_2 - 2 \tau_1) & 1 - \tau_1 & 3 - 2 \tau_1 \\ 
0 & 0 & \tau^2_2 (\tau_2 - \tau_1)^2 & p(\tau_1,\tau_2)& (1- \tau_1)^2 & 2 (1-\tau_1)(2-\tau_1)  \\ 
0 & 0 & 0 & \tau^2_2 (\tau_2 - \tau_1)^2 & (1-\tau_1)^2(1-\tau_2) & q(\tau_1,\tau_2)  \\
0 & 0 & (\tau_2 - \tau_1)^5 & 5 (\tau_2 - \tau_1)^4 & (1-\tau_1)^5 & 5 (1- \tau_1)^4  \\
0 & 0 & 0 & 0 & (1-\tau_2)^5 &   5 (1- \tau_2)^4  \\
\frac{\tau_1^6}{6!} 
& \frac{\tau_1^5}{5!} 
& \frac{\tau^6_2 - 2 (\tau_2 - \tau_1)^6}{6!} 
&  \frac{\tau_2^5 - 2 (\tau_2 - \tau_1)^5}{5!} 
& \frac{1 - 2(1-\tau_1)^6 + 2 (1-\tau_2)^6}{6!} 
& \frac{1 - 2(1-\tau_1)^5 + 2 (1-\tau_2)^5}{5!} 
%\end{array}
%\right]
%\end{eqnarray*}
\end{bmatrix} ,
$$
\normalsize
$p(\tau_1,\tau_2)= 2 \tau_2 (\tau_2 - \tau_1) (2\tau_2 - \tau_1)$, 
$q(\tau_1, \tau_2) =(1-\tau_1)(2(1-\tau_2)(2-\tau_1)+1-\tau_1)$.
Taking the $5$th derivative in \eqref{E3asDet}
and expanding along the last columns yields, for $0 \le t < \tau_1$,
$$
\mathcal E^{(5)}_3(t) = C (-5! \delta_1 + \delta_2 t),
$$
where $\delta_1$ and $\delta_2$ are the determinants of the submatrices of $D$ 
obtained by
removing the fourth row and the last row,  respectively.  
From $ \mathcal{E}_3^{(6)} (t) = -1$ for $0 \le t < \tau_1$, we derive $C = - 1/\delta_2$, and in turn  
\begin{eqnarray*}
\mathcal E^{(5)}_3(0) =  120 \,\frac{ \delta_1}{\delta_2}.
\end{eqnarray*}
In the case $\tau_2 = 2\tau_1$, an explicit calculation (facilitated by a computer algebra software) reveals that 
\begin{align*}
\delta_1 & = \frac{1}{360} \tau_1^{12} (1-2 \tau_1)^6 (4 - 32 \tau_1 + 189 \tau_1^2 -312 \tau_1^3 + 159 \tau_1^4),\\
\delta_2 & = 4 \tau_1^{13} (1-2 \tau_1)^6 (1-\tau_1) (7 - 5 \tau_1).
\end{align*}
Thus, as $\tau_1 \to 0$, we have
\begin{eqnarray*}
{\cal E}_3^{(5)} (0)  \sim  120 \frac{4 \tau_1^{12}/360}{28 \tau_1^{13} } 
= \frac{1}{21 \tau_1} \rightarrow + \infty.
\end{eqnarray*}
%\begin{eqnarray*}
%\mathcal E^{(5)}_3(0) = ? 
%\end{eqnarray*}
%In particular, by taking $\tau_2 = 2 \tau_1$ and letting $\tau_1 \searrow 0$ we get
%\begin{eqnarray*}
%\mathcal E^{(5)}_3(0) \sim  -\frac{1}{7 \times 360 \tau_1} \to -\infty 
%\end{eqnarray*}
This shows that Conjecture 2 does not hold.  \hfill $\Box$

\section{Alternative arguments}
\label{sec:AlternativeArguments}

Although the results of Section~\ref{sec:InterpErrorUnbounded} 
show that the methods of 
proof used in \cite{MR2382657}  (which are heavily based on the methods used in \cite{MR1041391} and \cite{MR1891742})
do not suffice for proving the desired rate results as stated there, we 
continue to believe that the rate will be $n^{-1/(2k+1)}$ for the ``gap conjecture'' of \cite{MR2382657}, and $n^{k/(2k+1)}$ for
the MLE of the $k-$monotone density $f_0$.  Here we sketch several possible routes toward proof of these conjectured results.

\subsection{Option A:  lower bound for the gaps}  
Note that the arguments in the preceding section showing unboundedness of the envelope of the 
interpolation error relied on taking $\tau_2 = 2\tau_1$ so that $\tau_2 - \tau_1 = \tau_1 \rightarrow 0$ 
where the $\tau$'s are regarded as parameters or variables indexing the entire class of interpolation errors 
for a scaling of the problem with 
$$
0 \equiv \tau_0 < \tau_1 < \cdots < \tau_{2k-4} < \tau_{2k-3} \equiv 1.
$$
Thus a ``coalesence'' of the knots leads to failure of the conjectures made in \cite{MR2382657}.

On the other hand, on the original time scale for the (random!) knots $\tau_0 < \tau_1 < \cdots < \tau_{2k-3}$ we want to show that
$\tau_{2k-3} - \tau_0 = O_p (n^{-1/(2k+1})$.  It seems likely that these random knots for the LSE actually do not ``coalesce'', but 
stay bounded away from each other asymptotically (at the rate $n^{-1/(2k+1)}$), and hence we expect to have
\begin{eqnarray}
\max_{1 \le j \le 2k+3} \frac{1}{ (\tau_j - \tau_{j-1})} = O_p (n^{1/(2k+1)} ) , 
\end{eqnarray}
or, equivalently
\begin{eqnarray}
\max_{1 \le j \le 2k+3} \frac{1}{n^{1/(2k+1)} (\tau_j - \tau_{j-1})} = O_p (1) .
\label{NonCoalesenceOfKnotsInProb}
\end{eqnarray}
If we could show that (\ref{NonCoalesenceOfKnotsInProb}) holds, then the classes of 
functions involved in the interpolation errors could be restricted to classes involving 
separated knots and the conjectures may be more plausible for these restricted classes.

\subsection{Option B:  alternative inequalities}  

While the methods of proof used in \cite{MR2382657} 
(and \cite{MR1041391}, \cite{MR1891742})
are based on empirical process inequalities which rely the small or scaling properties of 
envelopes (see e.g. Lemma A.1, page 2560 of \cite{MR2382657} or Lemma 4.1 of \cite{MR1041391}),
as opposed to smallness of the individual functions in the class relative to an envelope
as in Lemmas 3.4.2 and 3.4.3 of \cite{MR1385671} or \cite{MR2792551}.
On the other hand 
the proofs of the (global) rate of convergence of Hellinger distance 
from the  MLE $\widehat{f}_n$ to $f_0$ established in \cite{MR2520591} rely on 
inequalities for suprema of empirical processes based on uniform or bracketing entropy
for function classes in which the $L_2-$norms of individual functions are small relative to envelope 
functions (which may possibly be unbounded)  (see e.g. 
Lemmas 3.4.2 and 3.4.3 and Theorem 3.4.4 of \cite{MR1385671} for bracketing entropy 
type bounds, and see \cite{MR2792551} for classes with well behaved uniform entropy bounds; 
the results of \cite{MR2243881} might also be helpful in connection with the latter classes).

Thus there is some possibility that alternative inequalities for suprema of the empirical processes
involved may be needed in establishing the desired rate results when $k\ge 3$.  

\subsection{Option C:  alternative inequalities involving ``weak parameters''}  
While the inequalities discussed in option B above involve application of empirical process inequalities
involving ``strong parameters'' such as the expected values of envelope functions, there remains 
some possibility for the development of new inequalities based on ``weak parameters''; 
see e.g. the discussion on page 51 of \cite{MR2319879} 
and the material on page 209 of \cite{BoucheronLugosiMassart:13}.  
This option is the most speculative of the three.

%\nocite{*}
%\bibliographystyle{imsart-nameyear}
%\bibliographystyle{ims}
%\bibliography{SplineConjecture}

\end{document}